\documentclass[fleqn,10pt]{wlscirep}
\usepackage[utf8]{inputenc}
\usepackage[T1]{fontenc}
\usepackage{float}
\title{Identifying and Explaining the Resilience of Ecological Networks}

\author[1,*]{Cailan Jeynes-Smith}
\author[1,2]{Michael Bode}
\author[1,*]{Robyn P. Araujo}
\affil[1]{School of Mathematical Sciences, Queensland University of Technology, Brisbane, 4000,  Australia.}
\affil[2]{Safeguarding Antartica's Environmental Future, Queensland University of Technology, Brisbane, 4000, Australia.}

\affil[*]{cailan.jeynessmith@hdr.qut.edu.au, and r.araujo@qut.edu.au}


\begin{abstract}
Resilient ecological systems will be better able to maintain their structure and function in the emerging Anthropocene. Estimating the resilience of different systems will therefore provide valuable insight for conservation decision-makers, and is a priority goal of resilience theory. Current estimation methods rely on the accurate parameterisation of ecosystem models, or the identification of important motifs in the structure of the ecological system network. However, both of these methods face significant empirical and theoretical challenges. In this paper, we adapt tools developed for the analysis of biochemical regulatory networks to prove that a form of resilience - robust perfect adaptation - is a property of particular ecological networks, and to explain the specific process by which the ecosystem maintains its resilience. We undertake an exhaustive search for robust perfect adaptation across all possible three-species ecological networks, under a generalised Lotka-Volterra framework. From over 20,000 possible network structures, we identify 23 network structures that are capable of robust perfect adaptation. The resilient properties of these networks provide important insights into the potential mechanisms that could promote resilience in ecosystems, and suggest new avenues for measuring and understanding the property of ecological resilience in larger, more realistic socioecological networks. 
\end{abstract}

\begin{document}

\flushbottom
\maketitle

\thispagestyle{empty}

% \textbf{Short, Running Title:} Motifs of ecological resilience


% \textbf{Article Type:} Letter

\noindent \textbf{Author contributions statement}

CJS conceptualised the idea, RPA contributed to the modelling framework, CJS implemented methods, analysed results, created visualisations and produced the first draft of the manuscript, all authors contributed to the review and editing of the manuscript, and MB, and RPA supervised the project.

\noindent \textbf{Acknowledgements:} 

Computational resources and services used in this work were provided by the eResearch Office, Queensland University of Technology, Brisbane, Australia.

\noindent \textbf{Funding:} Robyn P. Araujo is the recipient of an Australian Research Council (ARC) Future Fellowship (project number FT190100645) funded by the Australian Government, Cailan Jeynes-Smith is supported by an Australian Government Research Training Program Scholarship.

\noindent \textbf{Conflict of Interest Statement:}

We declare we have no competing interests.

\noindent \textbf{Data Statement}

Data sharing is not applicable to this article as no new data were created or analyzed in this study. All code used to generate the results and figures in this article can be found on Github: \\
https://github.com/JeynesSmith/PerfectResilience.git

% \textbf{Number of Figures:} 5 \\
% \textbf{Number of References:} 76 \\
% \textbf{Abstract Length:} 197 words \\
% \textbf{Main Text Length:} 2454 words \\

\newpage 

\section*{Introduction}

\begin{figure}[H] % ht
\centering
\includegraphics[width=0.7\linewidth]{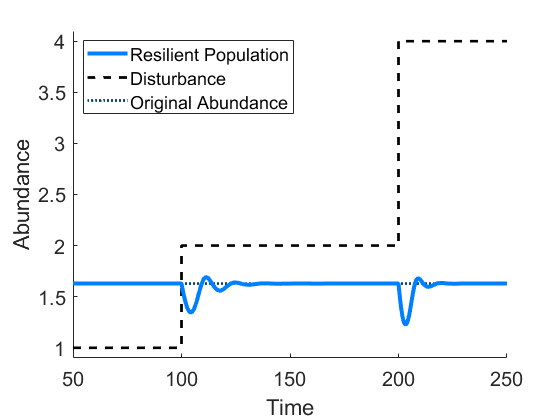}
\caption{An example of the perfect resilience behaviour, where a species' abundance is able to consistently return exactly to its original abundance following a disturbance. This behaviour is equivalent to robust perfect adaptation in biochemical reaction networks. Resilience, or imperfect adaptation are similar behaviours in which the abundance does not strictly return to the exact original abundance, but instead returns to within some surrounding region. This figure is generated by simulating a version of the network in Figure \ref{fig:Networks}(i). We simulate the system for one hundred time steps before doubling a stimulus ($S$) which disturbs the networks, specifically $S = \{1,2,4\}$. Parameters: $r_I=0.61$, $d_O=-0.07$, $a_{13}=-0.37$, $a_{21}=0.27$, $a_{22}=-0.81$, $a_{31}=0.32$, $a_{32}=-0.69$.}
\label{fig:ExampleAdapt}
\end{figure}

As the Anthropocene drives accelerating global change, resilience is an important and desirable characteristic of ecological and socioecological systems \cite{folke2004regime,may2019stability, parrott2012future}. Resilience has a range of definitions, but it primarily refers to the ability of an ecosystem to return to a particular set of ecological states following changes in environmental or social (policy) pressures (Figure \ref{fig:ExampleAdapt}). Resilience is generally used to describe a property exhibited by a particular ecological or socioecological system \cite{mumby2014ecological, urruty2016stability, evans2013robustness, digiano2012robustness, lake2013resistance, hoover2014resistance, macgillivray1995testing, holling1996engineering}, or as a measure of how well the population will recover from perturbation \cite{carpenter2001metaphor, lavorel1999ecological, holling1973resilience, neubert1997alternatives, meyer2016mathematical}. Understanding whether, and to what extent, an ecological or socioecological system is resilient helps us understand whether an ecosystem is at risk of collapse, or how far it can be pushed by environmental or social changes \cite{folke2004regime,may2019stability, parrott2012future}. 

A major hurdle to the use of resilience theory in conservation decision-making is its measurement. A primary approach is to create and parameterise a dynamical systems model of the ecological system, and then to simulate its response to perturbations \cite{mumby2014ecological, urruty2016stability}. However, the complexity of ecological systems, paired with sparse and noisy data sets, make the process of parameter estimation incredibly difficult \cite{adams2020informing, remien2021structural}. If resilience could be defined as a property of the structure of the ecological system, then parameter identifiability problems could be overcome. 

Previous studies into resilience of socioecological systems have proposed theoretical structures - network motifs - that underpin resilience \cite{barnes2017social, janssen2006toward}. However, this research is relatively atheoretical, with resilient network motifs being proposed on the basis of intuition, and justified by statistical association with observations of resilient dynamics. This limits the insights that these methods can offer.

By contrast, resilience has been studied extensively in biochemical reaction networks, where it is called adaptation. Resilience/adaptation is a relatively common dynamical property in cellular systems, where it has evolved to maintain the processes of life in a stochastic environment where external stimuli are constantly perturbing the abundances of molecules in the system. It has been observed in applications scaling from chemotaxis in single-celled organisms \cite{rao2004design, parent1999cell, macnab1972gradient, levchenko2002models, kollmann2005design, goy1977sensory, berg1972chemotaxis, alon1999robustness, yi2000robust} to complex sensory systems \cite{yau2009phototransduction, reisert2001response, matthews2003calcium, kaupp2010olfactory}, while the loss of adaptation has been linked to cancer progression, and substance abuse \cite{medina2009adaptation, araujo2007mathematical, fodale2011mechanism}. Resilience/adaptation is further categorised as either imperfect adaptation, where a particular element of the system returns within some accepted tolerance of its original (pre-stimulus) abundance \cite{ma2009defining,  frei2022genetic, aoki2019universal, briat2016antithetic, khammash2021perfect}, or as robust perfect adaptation, where that target element returns precisely to the exact same (pre-stimulus) abundance \cite{araujo2018topological, araujo2023design, araujo2023universal, jeynes2023protein, barkai1997robustness, sontag2003adaptation, ang2010considerations, ang2013physical, ferrell2016perfect, koshland1982amplification} (see Figure \ref{fig:ExampleAdapt}). Importantly, robust perfect adaptation is a property of the network interaction structure, as opposed to a property of a specific parameterisation of the network \cite{araujo2018topological}. 

In this study we apply analytical tools from biochemical reaction network theory to understand the process of resilience in ecological networks, and to assess whether ecological systems can exhibit robust perfect adaptation, and under what circumstances. To more closely match similar concepts in the existing ecological literature, we herein refer to `perfect resilience' as a mathematically equivalent property to robust perfect adaptation. We study ecological systems at the `operational layer' \cite{barnes2017social}, where the effect of policies directly impact populations in an ecosystem. We perform an extensive search of all three-species ecosystems using a novel application of the generalised Lotka-Volterra equations - a commonly used mechanistic framework in ecological modelling. We ask, are there network motifs that promote perfect resilience under an ecological framework, and will the corresponding networks still be representative of \textit{in situ} ecological systems? If possible, these specific motifs could be identified in ecosystems instead of current, more generalised motifs \cite{barnes2017social}.

\section*{Methods}

The structures capable of robust perfect adaptation (RPA) in biochemical reaction networks have recently been identified in full generality \cite{araujo2018topological, araujo2023universal}. Ma et al. \cite{ma2009defining} performed extensive numerical simulations of three-node chemical reaction (signalling) networks with Michaelis-Menten kinetics, and identified that there were two motifs structures that support RPA (at least approximately, as determined by tight numerical thresholds): a negative feedback loop with  buffer node (NFBLB), and an incoherent feedforward loop with proportioner node (IFFLP). Araujo and Liotta \cite{araujo2018topological} later determined that for networks of arbitrary size and complexity, and for arbitrary interaction kinetics, two-well defined subnetwork structures, or `modules', constitute a topological basis for RPA in any network:  Opposer modules, which are generalisations of three-node NFBLBs, and are feedback structured subnetworks;  and Balancer modules, which are generalisations of three-node IFFLPs, and are feedforward-structured subnetworks.  All RPA-capable networks, no matter how large or complex, and no matter the `kinetics' of the interacting elements, are necessarily decomposable into these two special modular subnetwork structures.  See \cite{araujo2018topological} for a detailed description of these mechanisms.  More recently, the intricate biochemical reaction structures that are compatible with these overarching RPA-conferring mechanisms have also been determined \cite{araujo2023universal}.

The biochemical networks that are capable of adaptation are often built from enzyme-mediated reactions \cite{ma2009defining, araujo2018topological, briat2016antithetic, aoki2019universal, jeynes2023protein}, where an enzyme combines with a substrate molecule to form an intermediate complex. This complex can either dissociate into the original molecules, or convert the substrate into a modified (product) form, releasing the enzyme unmodified. Importantly, there is no equivalent to these reactions in ecological systems, so it is unclear how, or if, ecological networks could generate adaptive behaviours without these fundamental reactions.

We provide an extensive study of all three-species interaction networks under a generalised Lotka-Volterra framework containing the three arbitrary species (or functional groups \cite{zheng1997soil}), $I$, $M$, and $O$, and a stimulus (cause of the change), $S$. We check for network structures in which the species, $O$, has the perfect resilience property. The stimulus represents an external influence on the system, which could be an intervention such as a harvest, or an environmental factor such as a heatwave, and can affect any combination of species in the system. 

\begin{figure}[H]
\centering
\includegraphics[width=0.8\linewidth]{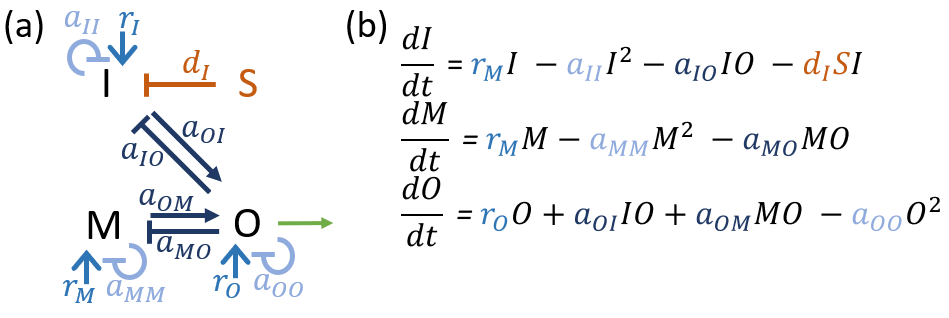}
\caption{(a) The graphical representation of a three-species network with the species, $I$, $M$, and $O$, and a stimulus $S$, and (b) the associated generalised Lotka-Volterra equations. This network depicts a system in which the species $O$ predates on both $I$ and $M$. Species $I$ is affected by the stimulus, and $O$ is the output of the network. In the network diagram (a) we represent a positive interaction by a pointed arrow and a negative interaction by a flat-ended arrow. The interactions and constants have been coloured based on the type of interaction: intraspecific competition, intrinsic growth rates, interspecific interactions (coloured light to dark blue respectively), and stimulus interactions (coloured red). The green arrow indicates that species $O$ is the output of the network, however this notation is dropped in later networks to simplify diagrams. }
\label{fig:GLVFramework}
\end{figure}

The full set of interactions in a three-species network modelled with the generalised Lotka-Volterra equations are given by,
\begin{align}
    \frac{\text{d}I}{\text{dt}} & = \overbrace{r_I I}^\text{Intrinsic Growth} - \overbrace{a_{II}I^2}^\text{Intraspecific Competition} + \overbrace{a_{IM} IM + a_{IO}IO}^\text{Interspecific Interactions} + \overbrace{d_ISI}^\text{Stimulus} , \label{GLVeqn1}\\
    \frac{\text{d}M}{\text{dt}} & = r_M M - a_{MM}M^2 + a_{MI} IM + a_{MO}MO + d_MSM, \\
    \frac{\text{d}O}{\text{dt}} & = r_O O - a_{OO}O^2 + a_{OI} IO + a_{OM}MO + d_OSO,\label{GLVeqn3}
\end{align}
where: $I$, $M$, and $O$ are the abundances of the interacting species (see Figure \ref{fig:GLVFramework}(a)); $r_i$ is the intrinsic growth rate of species $i$; $a_{ij}$ is the (\textit{per-capita}) interaction constant for how species $i$ is affected by species $j$; and $d_i$ is the interaction constant for how the stimulus, S, affects species $i$ ($d_i=0$ when species $i$ is not affected by the stimulus). 

The generalised Lotka-Volterra equations consist of three key terms: intrinsic growth, intra-specific competition (self-regulation), and inter-specific interactions. The intrinsic growth term models positive influences on a populations growth. This is often included for lower trophic-level species such as vegetation, where influences like rainfall and soil nutrients are implicitly modelled. Intra-specific competition represents a limiting or dampening term, whereby competition for the same resources within a population will ultimately limit its growth. Lastly, inter-specific interactions are all of the interactions between species $i$ and $j$. These can be positive or negative, where the signs of a pair of terms, $(a_{ij},a_{ji})$, are indicative of the type of interaction between species. Some common examples include: competition, $a_{ij},a_{ji}<0$; mutualism, $a_{ij},a_{ji}>0$; and predator-prey, $a_{ij}>0$, $a_{ji}<0$ where population $i$ is the predator. We illustrate an application of the generalised Lotka-Volterra equations and the graphical representation of a specific network in Figure \ref{fig:GLVFramework}.

We systematically examine the capacity for perfect resilience in networks constructed from every possible combination of terms in Equations (\ref{GLVeqn1})-(\ref{GLVeqn3}) based on a full factorial design. In the case where interactions are present or absent, i.e. not considering the sign of interactions, there are a total of $5\times 2^{12}=20,480$ possible network structures, since there are twelve parameters and five combinations in which the stimulus affects the network. When accounting for the sign of interactions, i.e. considering interactions as either absent, positive, or negative, there is a total of $20 \times 2^6\times 3^6 = 933120$ possible network structures. We therefore require efficient methods to test for the capacity of these networks to exhibit perfect resilience. 

Ma et al. \cite{ma2009defining} determined the capacity for imperfect adaptation in over 16,000 networks using metrics based on the simulated behaviour of 20,000 parameter sets for each network. This method is computationally expensive and depends heavily on an extensive sampling of parameter space for every individual network motif. In the current study, we automate algebraic methods recently developed by Araujo and Liotta \cite{araujo2023universal} which can definitively and accurately determine whether a given network motif has the capacity for RPA while handling all parameters symbolically. 

The cornerstone of this method is the recognition that any RPA-capable system, with dynamical rate equations $f_1, \ldots, f_n$ (assumed to be polynomial functions of the interacting elements, $x_1, \ldots, x_n$), is characterised by an {\it ideal} whose two-variable geometric projection assumes the special form of an RPA polynomial (see \cite{araujo2023universal} for full details).
In this case, RPA, and thus perfect resilience, requires the existence of three polynomials, $(p_1, p_2, p_3) \subset \mathbb{R}[I,M,O]$, such that
\begin{align}
    p_1 \frac{dI}{dt} + p_2 \frac{dM}{dt} + p_3 \frac{dO}{dt} = f(S,O)(O-k), \label{RPApoly}
\end{align}
where $k$ is a rational function of parameters, and $f(S,O)$ is a polynomial (known as the `pairing function' \cite{araujo2023universal}) that is non-vanishing on a suitably extensive region of the positive orthant.  The right-hand side of Equation (\ref{RPApoly}) has the form of an RPA polynomial.  If $(p_1,\ p_2,\ p_3)$ can be found that satisfy Equation (\ref{RPApoly}) for any given network, then the network in question has the capacity for perfect resilience, with the output steady-state value of $k$ for all disturbances and all parameter choices.  For the small network structures considered here, the existence of suitable $p_1,\ p_2,\ p_3$ can be determined algorithmically by computation of a Gr\"obner basis \cite{araujo2023universal} with an elimination monomial ordering, and with variables $S$ and $O$ ordered last.  We automate this process in Matlab using the \textit{lexicographic} monomial ordering, and then use Matlab's symbolic toolkit to automatically determine if a non-zero projection onto $S$ and $O$ exists, and if the projection can be factorised into an RPA polynomial (Equation (\ref{RPApoly})) based on the monomials. We reject any networks for which there is no non-trivial two-variable projection, or for which its two-variable projection is not an RPA polynomial.  All code developed for this study is provided at https://github.com/JeynesSmith/PerfectResilience.git.

We provide two examples in which a network does, and does not, have the capacity for perfect resilience based on the above projection test. The network in Figure \ref{fig:Networks}(i) (without $r_O$) is capable of achieving perfect resilience. By calculating the Gr\"obner basis for this network, we obtain the following projection,
$$ p_1 \frac{dI}{dt} + p_2 \frac{dM}{dt} + p_3 \frac{dO}{dt} = c_1 \mathbf{O^2S^2} + c_2 \mathbf{O^2S} + c_3 \mathbf{O^2} + c_4 \mathbf{OS^2} + c_5 \mathbf{OS} + c_6 \mathbf{O}, $$
where $c_i$ is some function of the network parameters. The projection can be factorised into a form which matches the right-hand-side of Equation (\ref{RPApoly}),
\begin{align}
    c_1 \mathbf{O^2S^2} + c_2 \mathbf{O^2S} + c_3 \mathbf{O^2} + c_4 \mathbf{OS^2} + c_5 \mathbf{OS} + c_6 \mathbf{O} = O(S^2+c_7S+c_8)(O-k), \label{CIC}
\end{align}
where $f(O,S)=O(S^2+c_7S+c_8)$ to match Equation (\ref{RPApoly}). Since the factorisation exists, this network has the capacity for perfect resilience. 

As a second example, the network from Figure \ref{fig:GLVFramework} has the following Gr\"obner basis projection,
\begin{align*}
    p_1 \frac{dI}{dt} & + p_2 \frac{dM}{dt} + p_3 \frac{dO}{dt} =  \\
    & c_1 \mathbf{O^5} + c_2 \mathbf{O^4S} + c_3 \mathbf{O^4} + c_4 \mathbf{O^3S^2} + c_5 \mathbf{O^3S} + c_6 \mathbf{O^3} + c_7 \mathbf{O^2S^2} + c_8 \mathbf{O^2S} + c_9 \mathbf{O^2} + c_{10} \mathbf{OS^2} + c_{11} \mathbf{OS} + c_{12} \mathbf{O} .
\end{align*}
There is no factorisation in which the right-hand-side of this equation will match the form of Equation (\ref{RPApoly}), even when based on the monomials ($S$ and $O$) alone, and therefore this network does not have the capacity for perfect resilience.

While the above condition is necessary for perfect resilience, we must still determine whether the setpoint, $O=k$, is a feasible and stable steady state. Feasibility ensures that all species have a positive abundance at steady state, while stability ensures that our ecosystem can move towards that steady state. Since the projection test has significantly reduced the number of networks (20480 to 1072 networks, see Supplementary Figure 3), we test for stability and feasibility using a sampling approach similar to Ma et al. \cite{ma2009defining} in which we randomly select $4 \times 10^3$ parameter sets. Parameters are selected from uniform distributions, $x \in (-1,1)$ for interspecific or stimulus interaction constants, $x \in (0,1)$ for intrinsic growth rates, and $x \in (-1,0)$ for intraspecific competition constants. We calculate the steady states of the network then substitute in random parameter values and check for feasibility i.e. there exists a steady state in which every species has a positive abundance. If the steady state is feasible, then we check stability using Lyaponuv stability i.e. negative real components for all eigenvalues of the Jacobian matrix \cite{lyapunov1992general}. If both stability and feasibility conditions are met, then the parameter set is saved as successful and we continue for the remaining $4\times 10^3$ parameter sets or until we obtain ten successful parameter sets. We define that a network is capable of perfect resilience if it has any successful parameter sets. Note, we aim to ensure that networks with perfect resilience do not require strict constraints on the parameters, and therefore do not require an intensive parameter search for stability and feasibility. 

We lastly ensure that successful systems have perfect resilience, and not a trivial form in which the output species has no reaction to changes in stimulus \cite{ma2009defining}. We simulate each system and check that the output reacts to a change in stimulus by at least $1\%$ of its pre-stimulus abundance and then returns to within $1\%$ of its pre-stimulus abundance at steady state. For networks which pass all of the above tests, we generate $2000$ parameter sets that enable feasible and stable steady states, and use these for further analysis. See Supplementary Figure 3 for a graphical representation of this process and the number of networks which proceed after each test.

\section*{Results}

In the following sections we identify all network configurations which have the perfect resilience property under the generalised Lotka-Volterra framework. We examine the structural constraints on ecosystem networks capable of perfect resilience, and how these constraints and the associated transient dynamics affect the possibility of observing perfect resilience in \textit{in situ} ecosystems.

\subsection*{Network Structures with Perfect Resilience} \label{Structure}

\begin{figure}[H]
\centering
\includegraphics[width=\linewidth]{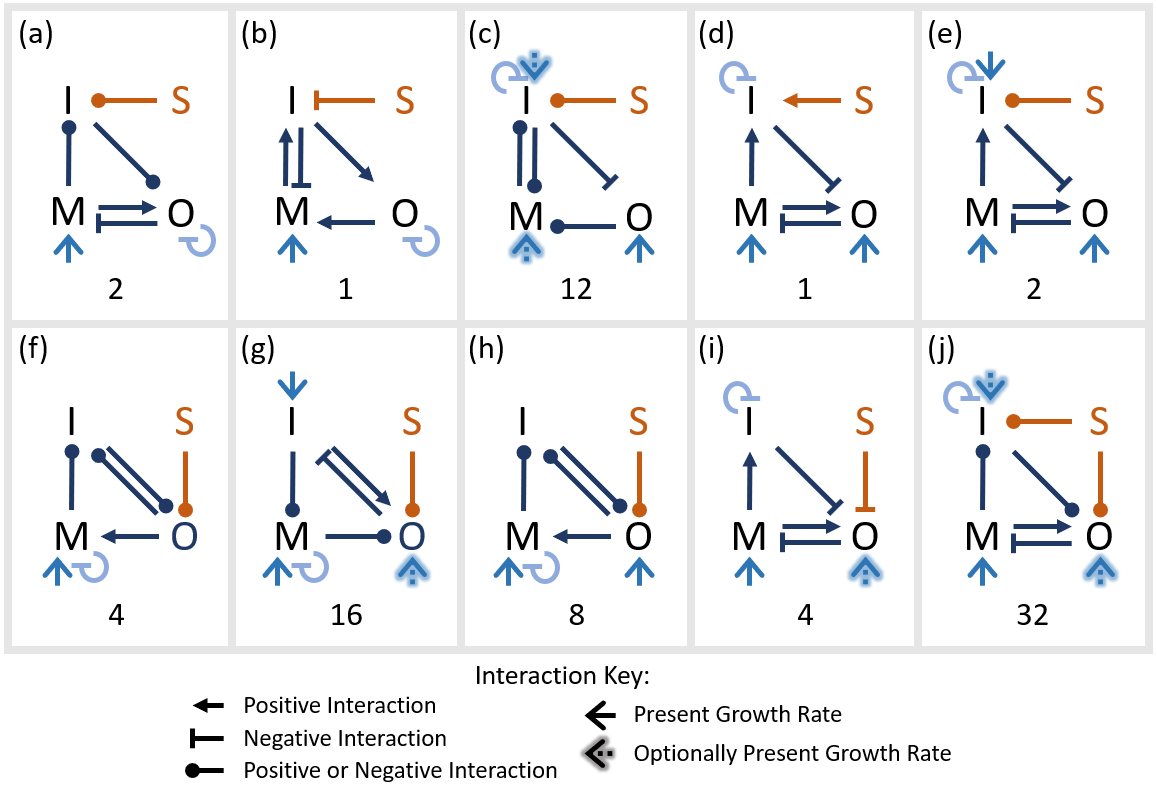}
\caption{Unique network motifs which enable perfect resilience in a generalised Lotka-Volterra three-species system. Below each motif is the number of variations on the motif which are capable of perfect resilience when the sign structure is accounted for - a total of 82 networks. The full list of 23 networks (without assigned sign structure) which are capable of perfect resilience can be found in Supplementary Figure 1. The 23 networks are condensed into ten motifs by including optionally present growth rates, where a growth rate may be present or absent and still have perfect resilience in the associated network. For example, (f) could be split into two networks with or without the $r_O$ growth rate. The motif in (c) is an exception where perfect resilience is present in networks in which: only $r_O$ is present; $r_O$ and $r_I$ are present; and when $r_O$, $r_I$ and $r_M$ are present. The motifs in (f)-(i) have non-unique variations in which species $I$ and $M$ are switched. The interaction signs are determined by 2000 parameter sets which permitted stable, feasible steady states, as discussed in the Methods section.}
\label{fig:Networks}
\end{figure}

Our extensive analysis of networks required testing a total of 20,480 possible network structures. In Supplementary Figure 1, we provide an extensive list of all 23 networks which are capable of perfect resilience before considering the sign of interactions between species. This is approximately $0.1\%$ of possible networks which are capable of perfect resilience. In Figure \ref{fig:Networks} we present the 23 networks as ten unique, condensed network motifs which illustrate the general trends in these structures.

\begin{figure}[H]
\centering
\includegraphics[width=0.9\linewidth]{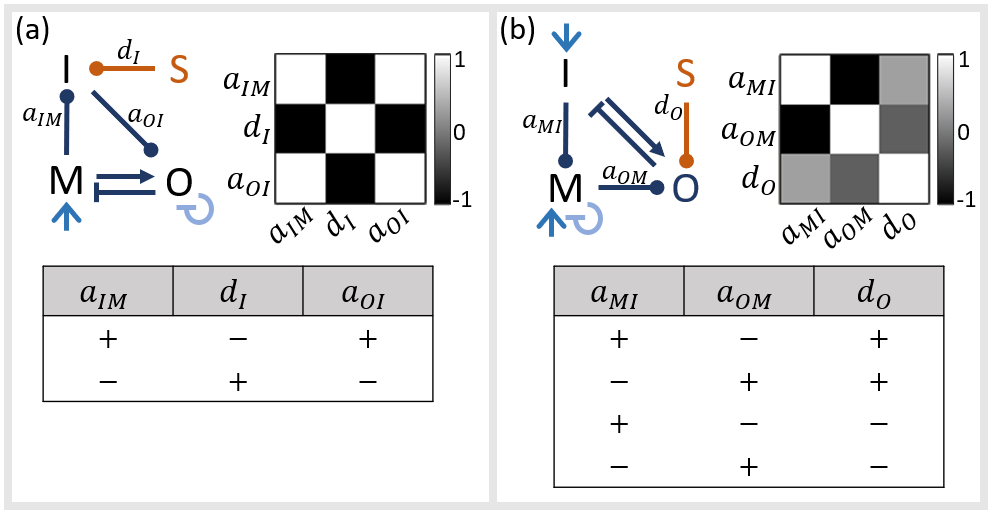}
\caption{Two examples of how the correlation between parameters can be used to specify the sign structure of networks capable of perfect resilience. The motifs in (a) and (b) are from Figure \ref{fig:Networks}(a) and (g) respectively. The correlation of parameter sets is represented in the grid where a value of $-1$ (black) represents a strong, negative correlation, and a value of $1$ (white) represents a strong, positive correlation. Strong correlations between parameters indicate a dependent relationship between the sign of those parameters, whereas a weak correlation (around $0$) indicates independence. These can then be used to create the associated tables of possible sign structure combinations.}
\label{fig:Correlation}
\end{figure}

After identifying the capacity for perfect resilience, we determine random parameter regimes which enable feasible and stable steady states for these networks (see Methods). Some parameters, like the growth rates and intraspecific competition constants, can only take on one sign (strictly positive or strictly negative respectively). To ensure feasibility and stability, some interaction and stimulus constants must take a specific sign, however most motifs had interactions that could take on positive or negative values (arrows with rounded heads in Figure \ref{fig:Networks}(a,c,e-h,j)). In these motifs, we studied the correlation between parameter sets (associated with feasibility and stability) and found that there is either a high correlation between the sign of interactions, or no correlation. This can then be used to further classify the motifs based on the sign of interactions - the sign structure. For example, the motif in Figure \ref{fig:Networks}(a) has three constants which can take on either sign: $a_{IM}$, $d_I$, and $a_{OI}$. These three constants have strong correlation which we use to classify the motif into two possible sign structures (Figure \ref{fig:Correlation}(a)). This can also be observed for the motif in Figure \ref{fig:Networks}(g), but one of the three interactions is weakly correlated with the others and therefore we can define four sign structures (Figure \ref{fig:Correlation}(b)). The number of sign structures are identified at the bottom of each motif in Figure \ref{fig:Networks}, and the specific sign structures can be found in Supplementary Figure 2. In total, when we specify the sign of every interaction, we identified 82 networks which had the perfect resilience property. These 82 networks represent less than $0.01\%$ of the possible network structures (when interaction sign is included).

While further classifying motifs based on sign structure can be used to provide more detailed motifs, we found that this did not provide distinguishable insights into the structure or dynamics of the ten unique network motifs in Figure \ref{fig:Networks} without the sign structure fully specified. In the following section, we draw conclusions on the structural requirements of these motifs and relate these to conditions we would expect from real ecosystems.

\subsection*{Topological Constraints and Ecosystem Implications}

In this extensive study, we examined the five possible configurations for how the stimulus can interact with populations in the network (excluding the sign and strength of those interactions), however only three were found to have motifs capable of perfect resilience when the stimulus affects: $I$, $O$, and $I$ and $O$ (Figure \ref{fig:Networks}, red arrows). Crucially there was no capacity for perfect resilience if the stimulus affects $I$ and $M$, or when it affects all three species. This is particularly important in the latter case as a stimulus representing environmental events, such as heatwaves or cyclones, are likely to directly affect all species in the network. This indicates that perfect resilience is never possible under these crucial types of stimuli which are increasing in frequency with climate change \cite{stott2016climate}.

When we translate the perfect resilience motifs into the mechanisms which obtain robust perfect adaptation in biochemical reactions, we identified that our motifs only used opposer mechanisms to generate perfect resilience - no balancer mechanisms. Araujo and Liotta \cite{araujo2018topological} determined that networks cannot have balancer mechanisms if the stimulus directly affects the output, because it requires multiple paths connecting the stimulus to the output to effectively `balance out' a change in stimulus. We do not observe multiple paths in any of the possible networks where the stimulus does not directly affect the output, $O$ (Figure \ref{fig:Networks}(a)-(e)). 

The structure of all ten motifs are dependent on a sequence of one-way interactions between species (Figure \ref{fig:Networks}, dark blue arrows). For example, the network in Figure \ref{fig:Networks}(a) has a one-way interaction connecting $M$ to $I$ followed by $I$ to $O$. In ecological systems having a one-way interaction is possible, an example being that an orchid on a tree benefits from the tree, but the tree has no significant benefit or harm from the orchid \cite{rasmussen2018epiphytic}. However, having \textbf{multiple} one-way interactions, in sequence, is highly unlikely in real ecological systems. 

Lastly, in all of the motifs we observe either intrinsic growth or self-regulation terms in only a subset of the species i.e. not for every species. It is not uncommon to exclude intrinsic growth terms for some species since, by definition, these terms are included to represent implicit increases for that population. However, self-regulation terms are more often included for all species as this ensures a carrying capacity for that species \cite{hening2018persistence, wangersky1978lotka, anisiu2014lotka}. In our perfect resilience motifs, only one species in the network ever has a self-regulation term (Figure \ref{fig:Networks}, light blue looping arrow), resulting in the potential for unbounded growth for two species. Self-regulation terms also play an important role in the stability of the steady state, and dampening oscillations of a species following a perturbation \cite{hening2018persistence, wangersky1978lotka, anisiu2014lotka}. In the following section we explore the dynamics of these networks and any implications that these have on perfect resilience.

\subsection*{Large Oscillations and Transient Dynamic Limitations}

In all ten motifs we observed a particular variety of integral control - constrained integral control \cite{araujo2023universal,xiao2018robust} - in which perfect resilience is present in the output species, $O$, as long as it does not go extinct. This occurs when there is an isolated factor of $O$ in the projection (Equation (\ref{RPApoly})), which can be observed in the example in Equation (\ref{CIC}). While this is realistic for an ecosystem, it is a stark contrast to chemical reaction networks where a molecule can be created from reactions independent of its abundance. Ecosystem networks are therefore limited by the strength of perturbations which the output is able to recover from while avoiding extinction.

\begin{figure}[H] % ht
\centering
\includegraphics[width=\linewidth]{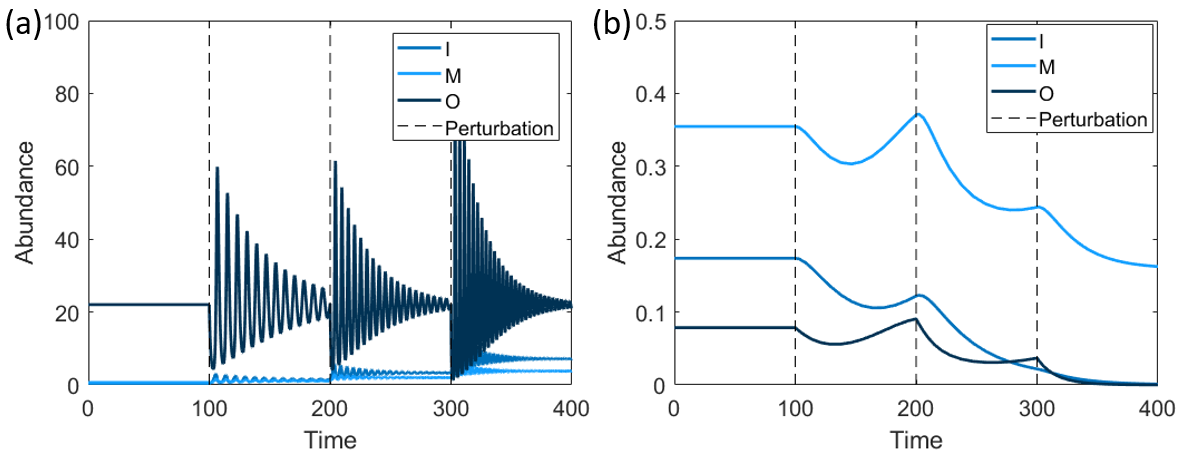}
\caption{Two examples of the oscillatory behaviours present in ecosystem networks capable of perfect resilience. In (a) we demonstrate a representative example of the highly oscillatory behaviour from the motif in Figure \ref{fig:Networks}(g) which is present in all motifs. Note that the oscillations have an increasing frequency and amplitude with repeated perturbations. In (b) we demonstrate how repeated perturbations before the network has returned to steady state results in extinction (motif from Figure \ref{fig:Networks}(h)). In both cases we simulate the system for one hundred time steps before doubling the stimulus, specifically $S = \{1,2,4,8\}$. Parameters: (a) $r_I=0.6$, $r_M=0.46$, $r_O=0.5$, $d_O=-0.74$, $a_{13}=-0.02$, $a_{21}=0.46$, $a_{22}=-1$, $a_{31}=0.85$, $a_{32}=-0.18$, (b) $r_M=0.04$, $r_O=0.09$, $d_O=-0.01$, $a_{12}=-0.18$, $a_{13}=0.83$, $a_{22}=-0.26$, $a_{23}=0.64$, $a_{31}=-0.42$.}
\label{fig:Oscillations}
\end{figure}

Generalised Lotka-Volterra models are prone to generating transient dynamics which are highly oscillatory \cite{wangersky1978lotka, anisiu2014lotka} and the networks which we identified as having perfect resilience are not exempt from this behaviour (Figure \ref{fig:Oscillations}(a)). When the stimulus changes, these networks rapidly oscillate in response and can take significant time to return to their pre-stimulus abundance. If another perturbation occurs within this oscillatory period it is possible that the abundance of the species can be perturbed to extinction and perfect resilience will be lost (Figure \ref{fig:Oscillations}(b)). Moreover, in reality having oscillations which repeatedly bring the output close to extinction have a higher risk of a stochastic event killing off the population. 

\section*{Discussion}

In this study we present a novel attempt to find perfect resilience in ecological systems. Moreover we have developed fully automated methods for identifying perfect resilience which are widely applicable to networks with any number of species, or functional groups, and under alternate modelling frameworks. We have identified that there is potential for ecological systems to obtain perfect resilience - albeit in only 23 network structures out of a possible 20,480 configurations in a three-species system. The motifs that we have identified can provide an important insight into ecological management decisions, particularly since this only requires an understanding of the structure of a network - an easier task compared to determining the interaction strengths in a network. By understanding the motif structure, it can then assist in highlighting interventions to ensure perfect resilience either does or does not occur for a target species or ecosystem function \cite{folke2004regime, may2019stability, parrott2012future}.

The 23 networks that we identified as being capable of perfect resilience were all based on opposer mechanisms \cite{araujo2018topological}, instead of the mixture of opposers and balancers which can promote robust perfect adaptation in biochemical reaction networks. This matches observations in biochemical literature, since the oscillatory behaviours of our networks are only possible under the feedback structures found in opposer mechanisms \cite{araujo2023design}. In biochemical networks there is a significantly better understanding of balancer mechanisms \cite{araujo2023design, goy1977sensory, mangan2003structure, goentoro2009incoherent, skataric2015fundamental, shinar2010structural}, compared to opposer mechanisms \cite{frei2022genetic, khammash2021perfect, briat2016antithetic, aoki2019universal, araujo2023design, jeynes2021ultrasensitivity, jeynes2023protein, ma2009defining, barkai1997robustness, koshland1982amplification, sontag2003adaptation, ang2010considerations, ang2013physical, drengstig2012basic, ferrell2016perfect}. There are two key drivers underlying opposer mechanisms - ultrasensitivity and antithetic integral control \cite{araujo2023design,  jeynes2023protein, ma2009defining, barkai1997robustness, sontag2003adaptation, ang2010considerations, ang2013physical, drengstig2012basic, ferrell2016perfect, koshland1982amplification, frei2022genetic, aoki2019universal, briat2016antithetic, khammash2021perfect}. These are dependent, respectively: on molecules rapidly converting to an active (`on') or inactive (`off') state, like a switch, to regulate an output; or the annihilation of excess output molecules, which get transformed into a combined (complexed) form before being removed. Neither of these drivers are examples of realistic ecosystem behaviours. Moreover, both drivers utilise the formation of enzyme-substrate complexes, but there is no equivalent to this set of reactions in species interactions. 

Many of the distinctions between ecological systems and biochemical reaction networks arise from structural differences between the two types of network. A clear feature in all ten network motifs, and likely pivotal to the opposer mechanism, is the presence of a sequence of one-way interactions between species. One-way interactions, where only one population is affected by an interaction, can be found in ecological networks, such as epiphytes living on trees \cite{rasmussen2018epiphytic}, but observing multiple one-way interactions in sequence is unlikely. Chemical reaction networks can include one-way interactions, dependent on the modelling framework used to describe the behaviour. For example, Michaelis-Menten equations readily display one-way interactions with enzymes which are able to alter molecular concentrations without having any effect on their own abundance \cite{ma2009defining}. In all of the perfect resilience motifs we also identified that only a subset of species have intraspecific competition terms. These terms play a crucial role in generalised Lotka-Volterra equations to dampen oscillations, and increase stability of the system \cite{hening2018persistence, wangersky1978lotka, anisiu2014lotka}. We observe that all of our motifs are associated with highly oscillatory behaviours which make them more susceptible to successive perturbations and at a higher risk of extinction from stochastic disturbances. Similar behaviours have been observed in coral reefs which are impacted by successive disturbances \cite{ortiz2018impaired}. Another crucial finding of this work is that there were no networks in which changes that affect all populations directly, can possibly have the perfect resilience property. This implies that disturbances, particularly environmental disturbances such as heatwaves, cannot be moderated by perfect resilience.

In this study, we have focused primarily on the commonly used generalised Lotka-Volterra equations and a stimulus which represents the direct effect of policy and environmental changes. An important problem with these equations are parameter identifiability \cite{adams2020informing, remien2021structural}. Generalised Lotka-Volterra equations simplify complex dimensions such as stochasticity, spatial and genetic heterogeneity, and complex interactions (e.g. prey switching). This makes it challenging to measure, or infer, the net \textit{per capita} strength of an interaction between any two species. While perfect resilience overcomes this problem by treating parameters as undefined variables, a more pressing concern is that there can be difficulties in even identifying the sign of the net interaction between species \cite{peterson2021reconstructing}. This poses a big problem with identifying these motifs in real ecosystems. While these problems are likely not limited to generalised Lotka-Volterra equations, it is important to note that the methods used here for identifying perfect resilience can readily be applied to other modelling frameworks, with only some constraints required for calculating the Gr\"obner basis \cite{araujo2023universal, buchberger2002grobner}.

The network motifs identified in this study have important implications for conservation management and socioecological systems. While these motifs can be used to identify resilient ecological networks, a pivotal problem that managers will face is attempting to design, or alter, ecosystems. In biochemical reactions identifying robust perfect adaptation can be used to suggest targeted treatments which can later undergo \textit{in vitro} or clinical testing. However, for managers, ecosystem engineering faces a slew of problems which cannot be tested before implementation. Removing species through culling can be controversial for native species \cite{read2021improving, drijfhout2022mind, drijfhout2020understanding, mckinnon2018media, mehmet2018kangaroo}, and the introduction of non-native species to an ecosystem has repeatedly had adverse effects \cite{fabres1978recent, nugent2001top, shine2020famous}. In this study, we only performed an extensive search of three-species ecosystems. An existing study from biochemical literature \cite{araujo2018topological} identifies how smaller mechanisms can be embedded within larger networks while ensuring that robust perfect adaptation, or perfect resilience, still obtains. There is still a benefit in searching for perfect resilience in larger networks to identify possible mechanisms that may only be possible with a larger number of species (potentially balancer mechanisms), and the methods developed here constitute an important foundation for further work in this area.  Future work may also be to identify the frequency with which these motifs occur \textit{in situ} \cite{barnes2017social}.

\bibliography{main}

\end{document}

\title{Identifying and Explaining the Resilience of Ecological Networks - Supplementary Materials}
\author{C. Jeynes-Smith, M. Bode, and R. P. Araujo}
\date{}
\maketitle

\subsection*{Supplementary Figure 1}

\begin{figure}[H] % ht
\centering
\includegraphics[width=\linewidth]{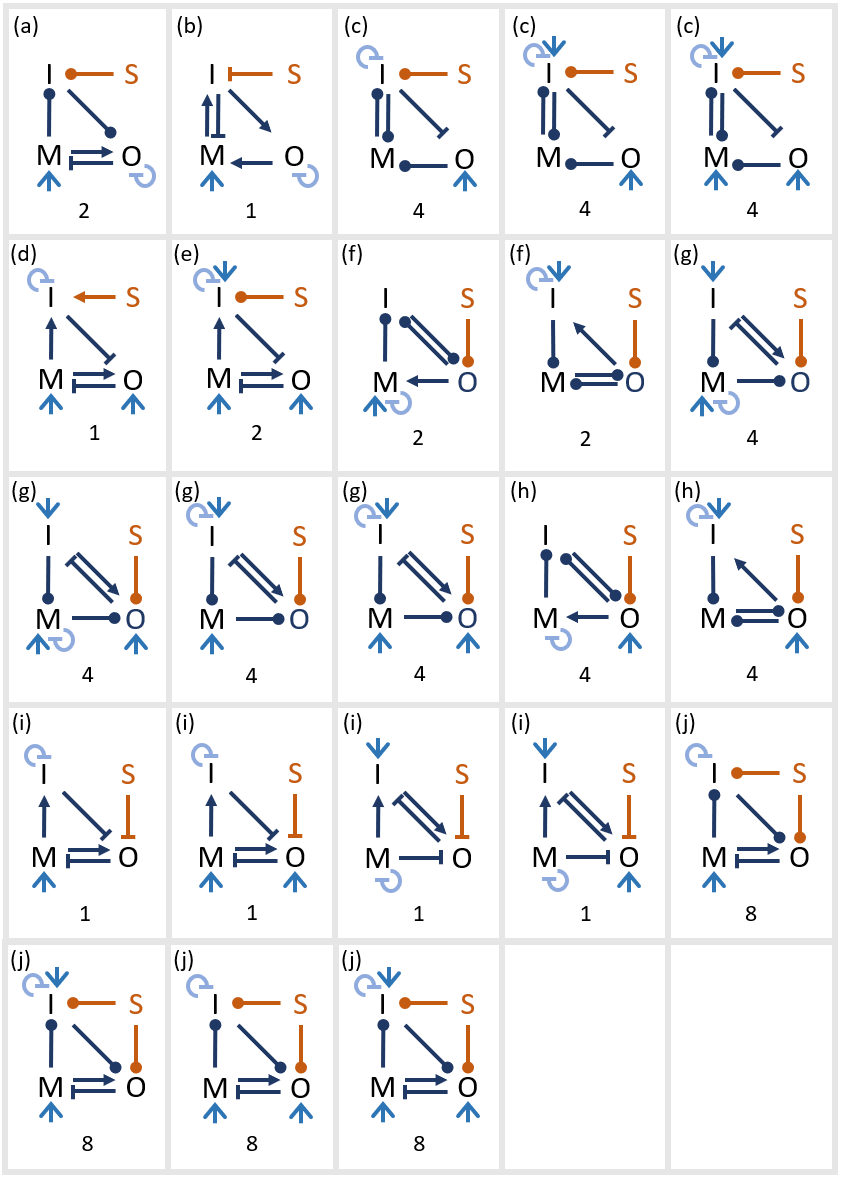}
\caption{Exhaustive list of all network structures capable of perfect resilience. This follows the same key as Figure 3 in the main text where positive interactions are a pointed arrow, negative interactions a flat arrow, and positive or negative interactions a ball-ended arrow. The number at the bottom of the network indicates the number of associated sign structures. Note that the extended forms in this figure (from the original motifs in Figure 3) have the same number of sign structures.}
%\label{fig:ExampleAdapt}
\end{figure}

\subsection*{Supplementary Figure 2}

\begin{figure}[H] % ht
\centering
\includegraphics[width=\linewidth]{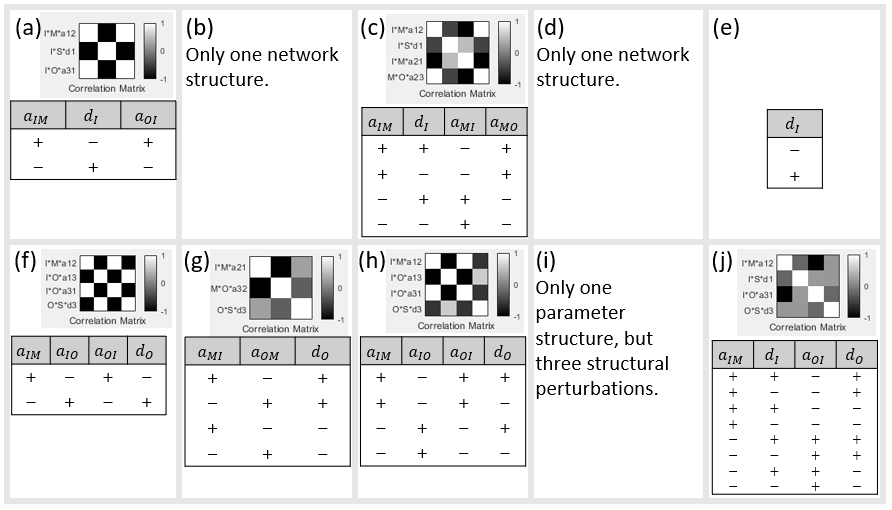}
\caption{Exhaustive list of correlation matrices for parameters with positive or negative signs from Figure 3. The correlation matrices are then used to define possible sign structures identified in the matching tables. Note that the alternate structures in Figure 3 which are obtained from added growth rates and switched species, all have the same correlation structure as the motif identified in Figure 3 and are therefore not separately listed.}
%\label{fig:ExampleAdapt}
\end{figure}

\subsection*{Supplementary Figure 3}

\begin{figure}[H] % ht
\centering
\includegraphics[width=\linewidth]{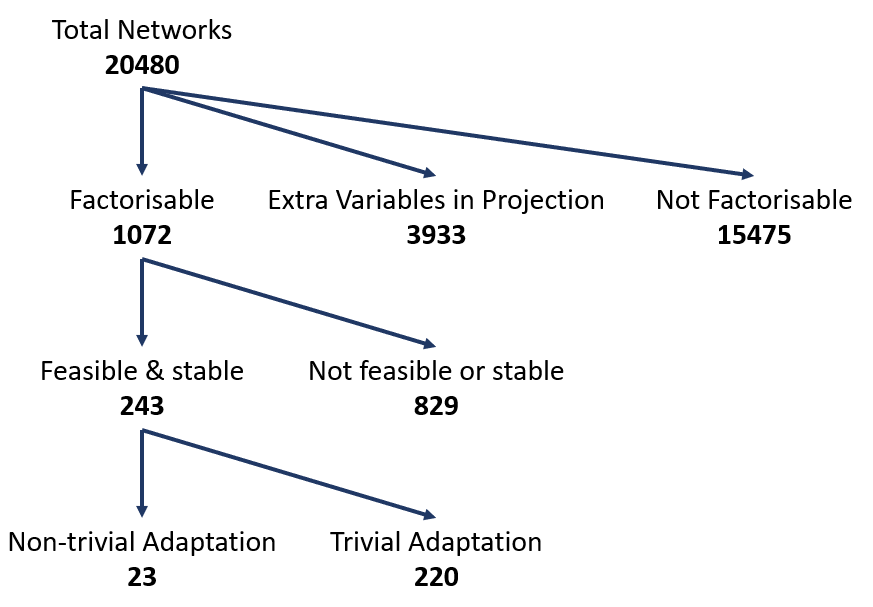}
\caption{A graphical representation of the methods used to identify perfect resilience in three-species Lotka-Volterra networks. The number of networks that fall under that definition is identified under the label. We initially start with 20480 networks, and end with only 23 networks that have perfect resilience.}
%\label{fig:ExampleAdapt}
\end{figure}